\documentclass{article}
\usepackage[english]{babel}
\usepackage{amssymb,amsmath,fullpage}
\usepackage{amsfonts,float,wrapfig}
\newtheorem{theorem}{Theorem}
\include{latexsymb}
\include{pxfonts}
\include{displaymath}
\include{tipa}
\begin{document}
\begin{center}
{\large \textbf{An alternative approach for\\compatibility of two discrete conditional distributions}}
\end{center}

\bigskip
\begin{center}
{\bf \large Indranil Ghosh}\\
{\large Austin Peay State University, USA}\\
e-mail: {\it  ghoshi@apsu.edu}\\
\mbox{}\\
\mbox{}\\
{\bf \large Saralees Nadarajah}\\
{\large  University of Manchester, United Kingdom}\\
e-mail: {\it mbbsssn2@manchester.ac.uk}\\
\end{center}

\bigskip
\begin{abstract}
Conditional specification of distributions is a developing area with increasing applications.
In the finite discrete case, a variety of compatible conditions can be derived.
In this paper, we propose an alternative approach to study the compatibility
of two conditional probability distributions under the finite discrete set up.
A technique based on rank based criterion is shown to be particularly convenient
for identifying compatible distributions corresponding to complete conditional specification including the case with zeros.
The proposed methods are illustrated with several examples.
\end{abstract}
\bigskip


\bigskip

\section{Introduction}

Specification of joint distributions by means of conditional densities has received considerable attention
in the last decade or so.
Possible applications may be found in the area of model building in classical
statistical settings and in the elicitation and construction of multiparameter prior
distributions in Bayesian scenarios.
For example, suppose that $\underline{Y} = \left(Y_{1},Y_{2}, \ldots, Y_{m}\right)$ is a $m$-dimensional random vector
taking on values in the finite range set
$\underline{\mathcal{Y}}_{1}\times \underline{\mathcal{Y}}_{2} \times \cdots \times \underline{\mathcal{Y}}_{m}$,
where $\underline{\mathcal{Y}}_{i}$ denotes possible values of $Y_{i}$, $i=1,2, \ldots, m$.
Efforts to ascertain
an appropriate distribution for $\underline{Y}$ frequently involve acceptance or rejection of a series of bets
about the stochastic behavior of $\underline{Y}$.
Suppose  we are facing a question of whether or not
to accept with odds five-to-one a bet that $Y_{1}$ is equal to one.
If we accept the bet then it puts a bound on the probability that $Y = 1$.
The basic problem is most easily visualized in the finite discrete case.

We contend that it is often easier to visualize that in some human population of interest,
the distribution of heights for a given weight will be unimodal with the mode of the conditional
distribution varying monotonically with weight.
Similarly, we may visualize a unimodal distribution of weights for a given height.
In that case, one might be interested to identify the joint nature of the height and weight data,
in the sense of deriving their joint probability distribution.
However, compatibility problems can be avoided or reduced if
one uses special types of models based on conditional probabilities.
Among them Bayesian network models (Lauritzen and Spiegelhalter, 1988; Pearl, 1988)
are strong candidates, due to their simplicity.
However, compatibility problem also arises when two or more experts participate in model selection processes,
or when we have partial information about conditional probabilities.

Several approaches exist in the literature with regard to the
problem of determination of the possible compatibility of two families of
discrete conditional distributions.
A review of the literature is provided in Section 2.

This paper focused on the finite discrete case takes a closer look at the compatibility problem,
viewing it as a problem involving linear equations in restricted domains.
Furthermore, we also focus our attention on situations where we have incomplete
(or partial) information on (either or both) of the two
conditional probability matrices, $A$ and $B$, under the compatible set-up.
In particular, we transform the problem of compatibility to a  linear programming (LP) problem
and we derive conditions for compatibility based on the rank of a matrix $D$, whose elements
are functions of the two conditional probability matrices,  $A$ and  $B$.
It will be shown that the problem of compatibility
is reduced to finding a solution to a set of $IJ$ equations in $(I-1)$ unknowns with
non negativity constraints, where $I$ and $J$ are the dimensions of the matrices $A$ and $B$.

This rank based criterion is different from Arnold et al. (1999) in which the authors
considered the rank of a matrix constructed in a different manner from the two conditional probability matrices $A$ and $B$.
Our method merits a separate study because we have transformed the problem of compatibility to a LP problem
and it's simplicity over other existing methods in the literature.
The purpose of this paper is to propose a general algorithm that can efficiently
verify compatibility in a straightforward fashion.
Our method is intuitively simple and general enough to deal with any full-conditional specification.

In situations where the two conditional probability matrices $A$ and $B$ have some
zeros appearing in the same position, one can not apply the cross product ratio criterion.
This is our main motivation to carry out this work.
We propose in Theorem 3 (Section 5) how  to verify the  existence of compatibility
in such situations and also to identify the existence of a unique joint probability matrix, $P$, based on the rank criterion.

We will mainly focus on the cases $I=2, 3$ and $J=2, 3$.
A detailed discussion will be provided for these cases and compatibility in higher dimensions will be discussed later on.
The paper proposes a convenient way to check whether or not two sets of conditional probabilities are
compatible based on the rank of a matrix $D$.
The paper also identifies all possible compatible distributions corresponding
to the given conditional probabilities in cases a unique solution fails to exist.

The contents of this paper are organized as follows.
In Section 3, we discuss  the concept of  compatibility
for any two conditional probability matrices $A$ and $B$.
In Section 4, we discuss in detail how the problem of
compatibility can be looked upon as a  LP problem.
In Section 5, we discuss the alternative approach to  compatibility
for the $(2\times 2)$ and $(3\times 3)$ cases.
In Section 6, some general discussion is made in the context of compatibility.
In Section 7, we discuss a real life example for illustrative purposes.

\section{A review}

Here, we review known work on compatibility of discrete conditional distributions.

Cacoullos and Papageorgiou (1983) provide characterizations of
discrete distributions by a conditional distribution and a regression function.
Arnold and Press (1989) derive if and only if conditions for existence of the joint distribution having given discrete conditionals.
Arnold and Gokhale (1994) characterize two-dimensional distributions representable as contingency tables in terms of uniform marginals.
Arnold and Gokhale (1998) use the Kullback-Leibler
information as an inconsistency measure for discrete
conditional distributions that are not compatible.
They provide algorithms computing the joint distribution that is
``least discrepant from the given inconsistent conditional specifications''.

Arnold et al. (2001) examine ``whether
a given set of constraints involving marginal and conditional probabilities
and expectations of functions are compatible''.
They also examine the case that the constraints are imprecisely specified.
Arnold et al. (2002) study when two families of discrete conditional distributions
are compatible, nearly compatible and $\epsilon$ compatible.
Arnold et al. (2004) present a method referred to as  ``rank one extension''
for deciding ``whether or not a set of conditional probabilities are compatible,
and if they are, obtaining all the associated compatible joint probabilities''.
The technique applies for both complete and partial conditional specifications.

Ip and Wang (2009) derive a canonical representation for multivariate distributions specified by discrete conditional distributions.
This representation gives if and only if conditions for existence and uniqueness of discrete conditional distributions.
Tian et al. (2009) present a unified method for checking
compatibility of discrete conditional distributions
that has connections to the ``quadratic optimization problem with unit cube constraints''.
Song et al. (2010) give if and only if conditions for compatibility of
discrete  conditional distributions and for uniqueness of the conditional distributions.
Wang and Kuo (2010) propose an algorithm for checking compatibility
of conditional discrete distributions having structural zeros.
They provide a Matlab program implementing their algorithm.
Kuo and Wang (2011) review the performance of five methods for checking
compatibility of discrete conditional distributions.
They also  propose an ``intuitively simple and general enough'' method for the same purpose.
Wang (2012) compares the performance of three methods for checking
compatibility of discrete conditional distributions.
Yao et al. (2014) propose a graphical approach for checking
compatibility of discrete conditional distributions.
They introduce a graphical representation
``where a vertex corresponds to a configuration
of the random vector and an edge connects two vertices if and only if the
ratio of the probabilities of the two corresponding configurations is specified
through one of the given full conditional distributions''.

\section{Compatibility}

Let $A$ and $B$  be two  $(I\times J)$ matrices with non-negative elements
such that $\displaystyle\sum_{i=1}^{I}a_{ij}=1$ for $j=1,\ldots,J$ and $\displaystyle\sum_{j=1}^{J} b_{ij} = 1$ for $i =1,2,\ldots,I$.
Without loss of generality, we assume that  $I\leq J$.
(If $I > J$ then replace $A$ and $B$ by the transpose of $A$ and the transpose of $B$, respectively.)
Matrices  $A$ and  $B$ are said to form a compatible conditional specification
(see Arnold et al. (1999))
for the distribution of $(X, Y)$ if there exists some $(I\times J)$ matrix $P$ with
non-negative entries  $p_{ij}$   and with $\displaystyle\sum_{i=1}^{I}\displaystyle\sum_{j=1}^{J}p_{ij}=1$ such that
\begin{eqnarray*}
\displaystyle
a_{ij}=\frac {p_{ij}}{p_{\cdot j}}
\end{eqnarray*}
and
\begin{eqnarray*}
\displaystyle
b_{ij}=\frac {p_{ij}}{p_{i \cdot}}
\end{eqnarray*}
for every $(i,j)$, where $p_{i \cdot} = \displaystyle\sum_{j=1}^{J}p_{ij}$ and $p_{\cdot j} = \displaystyle\sum_{i=1}^{I}p_{ij}$.
If such a matrix $P$ exists then, if we assume that
\begin{eqnarray*}
\displaystyle
p_{ij}=P\left(X=x_{i}, Y=y_{j}\right)
\end{eqnarray*}
for $i=1,2,\ldots, I$, $j=1,2,\ldots, J$, we have
\begin{eqnarray*}
\displaystyle
a_{ij}=P\left(X=x_{i}|Y=y_{j}\right)
\end{eqnarray*}
for $i=1,2,\ldots, I$, $j=1,2,\ldots, J$, and
\begin{eqnarray*}
\displaystyle
b_{ij}=P\left(Y=y_{j}|X=x_{i}\right)
\end{eqnarray*}
for $i=1,2,\ldots, I$, $j=1,2,\ldots, J$.

Equivalently,  $A$ and  $B$ are compatible if there exist stochastic vectors
$\underline{\tau} = \left( \tau_{1},\tau_{2},\ldots,\tau_{J} \right)$
and  $\underline{\eta} = \left(\eta_{1},\eta_{2},\ldots,\eta_{I}\right)$  such that
\begin{eqnarray*}
\displaystyle
a_{ij}\tau_{j}=b_{ij}\eta_{i}
\end{eqnarray*}
for every $(i, j)$.
In the case of compatibility, $\underline{\eta}$ and $\underline{\tau}$ can be readily interpreted
as the resulting  marginal distributions of $X$ and $Y$, respectively (see Arnold et al. (1999)).
For any probability vector $\eta = \left({\eta_1,\eta_2,\ldots, \eta_I}\right)$,  $p_{ij}=b_{ij}\eta_{i}$ is a probability
distribution on the $IJ$ cells.
So, the  conditional probability matrix denoted by $A$ has the elements
\begin{eqnarray}
\displaystyle
a_{ij}=\frac {p_{ij}}{\displaystyle\sum_{s=1}^{I}p_{sj}}=\frac {b_{ij}\eta_i}{\displaystyle\sum_{s=1}^{I}b_{sj}\eta_s}
\label{1}
\end{eqnarray}
for every $(i, j)$.

If $A$ and $B$ are compatible then
\begin{eqnarray*}
\displaystyle
a_{ij}\displaystyle\sum_{s=1}^{I}b_{sj} \eta_{s}=b_ {ij}\eta^{}_{i}.
\end{eqnarray*}
Then we have
\begin{eqnarray*}
\displaystyle
\tau_{j}=\displaystyle \sum_{s=1}^{I}b_{ij}\eta_{s}
\end{eqnarray*}
for $j=1,\ldots, J$.
In that case, expressions given in (\ref{1}) can be written as
\begin{eqnarray*}
\displaystyle
a_{ij}\displaystyle\sum_{s=1}^{I}b_{sj}\eta^{}_{s}-b_{ij}\eta^{}_{i}=0.
\end{eqnarray*}

In matrix notation the above can be written as
\begin{eqnarray}
\displaystyle
D\underline{\eta}=0,
\label{DD}
\end{eqnarray}
where $D$ is a matrix of dimension $IJ\times I $ and
$D\underline{\eta}=0$ is a  homogeneous system of  $IJ$  equations in $I$ unknowns  $\eta_{i}$.
Through well known matrix operations (such as left-multiplication by non-singular matrices)
it's rows can be reduced to at most  $I$ rows with non-zero elements (the so called ``Row Echelon form'').
Now let this reduced system be denoted by  $D_{r}\underline{y} = 0$,
where $\underline{y} = \left(y_{1},y_{2},\ldots, y_{I}\right)'$.
Matrices $A$ and $B$ are compatible if the system $D_{r}y = 0$ has a solution $y^{*}$ of non-negative elements with at
least one positive element.
If such a $y^{*}$ exists it can be scaled to arrive at a probability vector.
However, $A$ and $B$ are not compatible if the only solution with non-negative elements of  $D_{r}y=0$ is the null vector.
In order to examine whether or not such a solution $y^{*}$ of  $D_{r}y = 0$  exists (especially when $I$ is large),
one could solve the feasible region alone or the methodology of LP may be used.
Specifically, consider the problem of
maximizing the objective function $\displaystyle\sum_{i}y_{i}$ subject to
(a) the non-negativity constraints $\displaystyle\sum_{i=1}^{I}y_{i}\geq 0$,
(b) the equality constraints  $D_{r}y=0$,
and (c) the constraint $\displaystyle\sum_{i}y_{i=1}^{I}\leq 1$.
If the maximum of the objective function is positive, then the corresponding
optimizing vector is $y^{*}$, which can be scaled into the probability vector $\eta^{*}$.
If the maximum is zero, then $A$ and $B$ are not compatible.
We now consider some results related to
compatibility of two matrices $A$ and $B$ and later on
we discuss an alternative approach to  compatibility.

\section{Compatibility of two matrices $A$ and $B$}

We know that if the matrices $A$ and $B$ are compatible then $a_{ij}p_{\cdot j}=b_{ij}p_{i \cdot}$
for every $i = 1, 2, \ldots, I, j = 1, 2, \ldots, J$
(see Arnold et al. (1999)).
Equivalently, we can write
\begin{eqnarray*}
\displaystyle
a_{ij}\displaystyle\sum_{s=1}^{I}p_{sj}-b_{ij}\displaystyle\sum_{k=1}^{J}p_{ik}=0
\end{eqnarray*}
for every $i = 1, 2, \ldots, I, j = 1, 2, \ldots, J$,
which again can be written as
\begin{eqnarray*}
\displaystyle
a_{ij}\left[p_{1j}+p_{2j}+\cdots+p_{ij}+\cdots+p_{Ij}\right] -
b_{ij}\left[p_{i1}+p_{i2}+\cdots+p_{ij}+\cdots+p_{iJ}\right]=0
\end{eqnarray*}
for every $i = 1, 2, \ldots, I, j = 1, 2, \ldots, J$.

In matrix notation, the above system of linear equations can be written as
\begin{eqnarray}
\displaystyle
C\underline{p}=0,
\label{111}
\end{eqnarray}
where $C$ contains elements calculated from those of $A$ and $B$ and is a matrix of dimension $IJ\times IJ$ and
$\underline{p}^{(IJ\times 1)} = \left(p_{11},p_{12},\ldots,p_{IJ}\right)^{T}$.
Compare (\ref{111}) to (\ref{1}).
We have the following theorem.
\bigskip

\begin{theorem}
The solution space,  $\Omega$, for the system of equations in (\ref{1})
is $(I-M)\underline{z}$,  where $M$ is an idempotent matrix and $\underline{z}^{(IJ\times 1)}$
is any arbitrary vector of dimension $IJ\times 1$.
\end{theorem}

\noindent
{\bf Proof:}
We may consider  $M=C^{-}C$,  where  $C^{-}$ is the $g$-inverse of the matrix $C$.
We have considered the
$g$-inverse of  $C$  since $\text{\rm rank} \left(C^{IJ\times IJ}\right) < IJ$.
Next observe that
\begin{eqnarray*}
\displaystyle
M
&=&
\displaystyle
C^{-}C,
\nonumber
\\
\displaystyle
M^{2}
&=&
\displaystyle
C^{-}CC^{-}C
\nonumber
\\
&=&
\displaystyle
C^{-}C
\nonumber
\\
&=&
\displaystyle
M,
\end{eqnarray*}
which follows from the definition of $g$-inverse  since $CM=CC^{-}C=C$.
Hence, each of the $IJ$ columns
$\underline{h}_{1}, \underline{h}_{2},\ldots, \underline{h}_{IJ}$
of  $(I-M)$  is orthogonal to the rows of $C$.
But
\begin{eqnarray*}
\displaystyle
(I-M)^{2}=I-M-M+M^{2}=I-M
\end{eqnarray*}
since $M^{2}=M$.
Also
\begin{eqnarray*}
\displaystyle
\text{\rm rank}(I-M)=\text{\rm tr}(I-M)=\text{\rm tr}(I)-\text{\rm tr}(M)=IJ-r,
\end{eqnarray*}
where $r=\text{\rm rank}(C) = \text{\rm rank(M)}$ and tr denotes trace.
So, only $(IJ-r)$ of the column vectors
$\underline{h}_{1}, \underline{h}_{2}, \ldots, \underline{h}_{IJ}$
are linearly independent and without loss of generality we may consider them to be
$\underline{h}_{1},\underline{h}_{2},\ldots,\underline{h}_{IJ-r}$.

Again since  $C$  is an  $(IJ\times IJ)$  matrix of rank $r$, it's rows are $IJ$-vectors
and therefore we can find at most  $(IJ-r)$ linearly independent vectors orthogonal to them and
$\underline{h}_{1},\underline{h}_{2},\ldots,\underline{h}_{IJ-r}$
is one such set.
If there is any other vector orthogonal to the rows of  $C$,  it must be a linear combination
of  $\underline{h}_{1},\underline{h}_{2},\ldots,\underline{h}_{IJ-r}$.
Now since  $C\underline{p}=0$,  any vector  $\underline{p}$  satisfying  $C\underline{p}=0$
must be a linear combination
of $\underline{h}_{1},\underline{h}_{2},\ldots,\underline{h}_{IJ-r}$.

But equivalently we can say that  $\underline{p}$  is a  combination of
$\underline{h}_{1}, \underline{h}_{2}, \ldots, \underline{h}_{IJ}$ because
$\underline{h}_{IJ-r+1}, \underline{h}_{IJ-r+2}, \ldots, \underline{h}_{IJ}$
are  combinations of  $\underline{h}_{1}, \underline{h}_{2}, \ldots, \underline{h}_{IJ-r}$.
Hence, $\underline{p}$ must be of the form
\begin{eqnarray*}
\displaystyle
\underline{p}
&=&
\displaystyle
z_{1}\underline{h}_{1}+z_{2}\underline{h}_{2}+z_{3}\underline{h}_{3}+\ldots+z_{IJ}\underline{h}_{IJ}
\nonumber
\\
&=&
\displaystyle
\left(\underline{h}_{1},\ldots,\underline{h}_{IJ}\right)z
\nonumber
\\
&=&
\displaystyle
(I-M)\underline{z}
\end{eqnarray*}
for some choices of  $\underline{z} = \left(z_{1}, \ldots, z_{IJ}\right)'$.

Conversely, if the above holds then
\begin{eqnarray*}
\displaystyle
CM
&=&
\displaystyle
C(I-M)\underline{z}
\nonumber
\\
&=&
\displaystyle
(C-CM)\underline{z}
\nonumber
\\
&=&
\displaystyle
0
\end{eqnarray*}
because of the fact that  $CM=C$.
Hence, the proof.
$\square$

Next we have the following two propositions:
\begin{itemize}

\item
\noindent
{\bf Proposition 1:}
If there exists a vector in $\Omega$ such that
all its elements are non-negative and at least one element is
strictly positive then  $A$  and  $B$ are compatible.

\noindent
{\bf Proof:}
From our earlier result, the solution space is
\begin{eqnarray*}
\displaystyle
\Omega=(I-M)\underline{z},
\end{eqnarray*}
where $M=C^{-}C$.
Now suppose  $\underline{z} = \left(z_{1}, z_{2}, \ldots, z_{J}\right)'$,   where  $z_{u} = \left(z_{1u}, \ldots, z_{Iu}\right)'$.
Let us consider  $z_{i}\geq 0$ and (say)  $z_{j}>0$, where  $j\in i$,
meaning at least one of the  $z_{i}$'s  is strictly zero.
Then we can write  $W = \displaystyle\sum_{u=1}^{J}z_{Iu}$.
Obviously, $W>0$.
So, we can rewrite $\underline{z}$ as
\begin{eqnarray*}
\displaystyle
\underline{z}_{\rm new}
&=&
\displaystyle
\left(\frac {z_{1}}{W},\frac {z_{2}}{W},\ldots,\frac {z_{J}}{W}\right)
\nonumber
\\
&=&
\displaystyle
\left(p_{1},p_{2},\ldots,p_{J}\right)',
\end{eqnarray*}
which is a valid probability distribution since $\displaystyle\sum_{i=1}^{I}\displaystyle\sum_{j=1}^{J}\frac {z_{ij}}{W}=1$.
Hence, the proof.
$\square$

\item
\noindent
{\bf Proposition 2:}
If every vector in $\Omega$ has non-positive elements then  $A$  and  $B$ are not compatible.

\noindent
{\bf Proof:}
In this situation suppose that
$\underline{z}^{(IJ\times1)} = \left( \underline{z}_{1}, \underline{z}_{2}, \ldots, \underline{z}_{J} \right)'$
with  $\underline{z}_{i}\geq 0$ for $i = 1, 2, \ldots, J$.
But if we consider the trivial situation or the possibility of every $\underline{z}_{i}=0$ then all the elements of
$\underline{z}^{(IJ\times 1)}$ are zero which implies that it can not be a valid probability distribution.
$\square$
\end{itemize}

Next observe that since $\displaystyle\sum\displaystyle\sum p_{ij}$ is a linear function in $p_{ij}$,
the problem is equivalent to the following LP problem:
Maximize  $f \left(\underline{p}\right) = \displaystyle\sum\displaystyle\sum p_{ij}$
subject to $C\underline{p}=0$ and $1\geq p_{ij}\geq 0$ for every $(i,j)$.

\begin{theorem}
In the above LP problem, max $f \left(\underline{p}\right) > 0$ if and only if $A$ and $B$ are compatible.
\end{theorem}

\noindent
{\bf Proof:}
Note that if max $f \left(\underline{p}\right) > 0$ then at least one of the elements
in $\underline{p} = \left(p_{11},\ldots,p_{IJ}\right)'$ is strictly positive,
so we may consider  $p_{uv}>0$, where $u\in i, v\in j$ and $p_{ij}\geq 0$ for $u\neq i$, $v\neq j$.
Then obviously we have
\begin{eqnarray*}
\displaystyle
{\rm max} f\left(\underline{p}\right)=
\displaystyle\sum_{i}
\displaystyle
\sum_{j}p_{ij}=
\displaystyle\sum_{i\neq u}
\displaystyle
\sum_{j\neq v}p_{ij}+p_{uv}>0.
\end{eqnarray*}
Next note that we can write
$\underline{p} = \Big(\underbrace{p_{11},p_{12},\ldots,p_{1J}}$,
$\underbrace{p_{21},p_{22},\ldots,p_{2J}},\ldots,\underbrace{p_{I1},p_{I2},\ldots,p_{IJ}} \Big)$.
Then we can rewrite $\underline{p}$ as
\begin{eqnarray*}
\displaystyle
\underline{p}
&=&
\displaystyle
\left(\underline{p}_{1},\underline{p}_{2},\ldots,\underline{p}_{I}\right)'
\nonumber
\\
&=&
\displaystyle
\left(\frac {\underline{p}_{1}}{{\rm max} f(p)},
\frac {\underline{p}_{2}}{{\rm max} f(p)},
\ldots,
\frac {\underline{p}_{I}}{{\rm max} f(p)}\right).
\end{eqnarray*}
Then it  becomes a valid probability distribution.
We have from our previous
result that $A$ and $B$ are compatible.
Hence, the proof.
$\square$

\section{An alternative approach to  compatibility}

Questions  of compatibility of conditional and marginal specifications of distributions are of fundamental
importance in modeling scenarios.
May be the earliest work in this area is Patil (1965).
He considered the discrete case under a mild regularity condition and showed that the conditional distributions
of $X$ given $Y$ and $Y$ given $X$ will uniquely determine the joint distribution  of $(X, Y)$.
There are several versions of necessary and sufficient conditions for compatibility given
by Arnold and Press (1989) and Arnold et al. (2002, 2004).
In some situations, the Arnold et al. (2004) condition for checking compatibility
was found to be difficult and less effective.
This is why we came up with a relatively easy and simple procedure to check the condition for compatibility.
The new method, which requires only some elementary type operation of matrices (``Row Echelon form'')
could provide a much simpler and a more effective approach.

When the given conditional distributions are compatible,
it is natural to ask whether the associated joint distribution is unique.
This issue has been addressed in the literature by Amemiya (1975),
Gourieroux and Montfort (1979), Nerlove and Press (1986) and Arnold and Press (1989).
Arnold and Press (1989) pointed out that the condition for uniqueness is generally
difficult to check.
In this paper, through the structure of the reduced $D$ matrix,
we provide a simple criteria for checking uniqueness.

Here, we discuss the compatibility of two conditional matrices $A$ and $B$ along with the uniqueness
and the existence of a joint probability $P$ based on the rank of a matrix $D$.
The salient feature of Theorem 3 lies in the fact that it can be applied to situations
where matrices $A$ and $B$ have some zeros appearing in the same position.
In situations like this, the cross product criterion can not be applied to check compatibility.
This is the major motivating factor for this paper.

Now, we consider the theorem:

\begin{theorem}
Any two given conditional probability matrices $A$ and $B$ of dimension $(I\times J)$
are compatible if $\text{\rm rank} \left(D^{(IJ\times I)} \right) \leq I-1$ with equality when there
exists a unique solution for the unknown $\eta_{i}$ for every $i$.
\end{theorem}

\noindent
{\bf Proof:}
Note that $\text{\rm rank}\left(D^{(IJ\times I)}\right) \leq \text{\rm min}(IJ,I)=I$.
Now when $D$ has full rank, i.e., $\text{\rm rank}(D)=I$, the only
solution to $D\underline{\eta}=0$ is the null vector (trivial solution).
So, matrices $A$ and $B$ are incompatible.

Next if we have  $\text{\rm rank} \left( D^{(IJ\times I)} \right) \leq I-1$, the number
of equations $(IJ)$ is greater than number of unknowns $(I-1)$, so we must have a non-trivial solution.
If the non-trivial solution is positive then it can be appropriately scaled to arrive
at a probability vector $\underline{\eta}^{*}$.
Hence, the two matrices $A$ and $B$ are compatible.
However, in this case, the system of equations is not homogeneous and we have at most $(I-1)$ solutions.

When  $\text{\rm rank}(D)=I-1$, we have $(I-1)$ unknowns
subject to the linear constraint $\displaystyle\sum_{i=1}^{I}\eta_{i}=1$.
The $(I-1)$ equations
(excluding the redundant equations from the total set of $IJ$ equations) and the system
of linear equations is homogeneous so that there exists a unique solution.
This completes the proof.
$\square$

This theorem is useful in situations when the two conditional matrices $A$ and $B$ have
zeros as elements appearing in the same position and we can not guarantee the existence of a compatible
matrix $P$ by the cross product ratio criterion.

Next, we discuss  compatibility for  $(2\times 2)$, $(3\times 3)$  cases with examples.
In the examples of Sections 5.2 and 5.3, ranks of all numerical matrices were determined
by using the Mathematica software.
The construction of the matrix $D$ via (\ref{DD}) and finding its rank via Mathematica is easy.
This is especially so since Mathematica is more powerful than
other software like Maple and Matlab and is widely available
in almost every mathematics/statistics department.

\subsection{Compatibility in $(2\times 2)$ case}

\noindent
{\bf Corollary 1:}
If the two matrices $A$ and $B$   are of  dimension $2\times 2$,
then rank $(D) = 1$, in case $A$ and $B$ are compatible.

\noindent
{\bf Proof:}
In a $(2\times 2)$ case, the $D$ matrix is
\begin{eqnarray*}
\displaystyle
D=\left(
\begin{array}{cc}
\displaystyle
b_{11}\left(a_{11}-1\right) & a_{11}b_{21}\\
\displaystyle
b_{12}\left(a_{12}-1\right) & a_{12}b_{22}\\
\displaystyle
a_{21}b_{11}& b_{21}\left(a_{21}-1\right) \\
\displaystyle
a_{22}b_{12}  &b_{22}\left(a_{22}-1\right) \\
\end{array}
\right).
\end{eqnarray*}
Equivalently, we can write $D$ as (because $a_{11}+a_{21}=1$ and $a_{12}+a_{22}=1$)
\begin{eqnarray*}
\displaystyle
D=\left(
\begin{array}{cc}
\displaystyle
-b_{11}\left(a_{21}\right) & a_{11}b_{21}\\
\displaystyle
-b_{12}\left(a_{22}\right) & a_{12}b_{22}\\
\displaystyle
a_{21}b_{11}  & -b_{21}\left(a_{11}\right) \\
\displaystyle
a_{22}b_{12}  & -b_{22}\left(a_{12}\right) \\
\end{array}
\right).
\end{eqnarray*}
Next we consider the elementary row-transformations:
\begin{itemize}

\item
new(row 3)=row 1 + row 3;

\item
new(row 4)=row 4 + row 2.

\end{itemize}
So, our $D$ matrix reduces to
\begin{eqnarray*}
\displaystyle
D=\left(
\begin{array}{cc}
\displaystyle
-b_{11}\left(a_{21}\right) & a_{11}b_{21}\\
\displaystyle
-b_{12}\left(a_{22}\right) & a_{12}b_{22}\\
\displaystyle
0  & 0 \\
\displaystyle
0  & 0 \\
\end{array}
\right).
\end{eqnarray*}
However, in a $(2\times 2)$ case, if $A$ and $B$ are compatible, we can write
\begin{eqnarray}
\displaystyle
a_{12}a_{21}b_{22}b_{11}=a_{11}a_{22}b_{21}b_{12}.
\label{3}
\end{eqnarray}
Again we apply the elementary row transformations:
\begin{itemize}

\item
new(row 1)=row 1$\times \left(a_{12}b_{22}\right)$;

\item
new(row 2)=row 2$\times \left(a_{11}b_{21}\right)$.

\end{itemize}
So, our $D$ matrix reduces to
\begin{eqnarray*}
\displaystyle
D=\left(
\begin{array}{cc}
\displaystyle
-b_{11}\left(a_{21}\right)a_{12}b_{22} & a_{11}b_{21}a_{12}b_{22}\\
\displaystyle
-b_{12}\left(a_{22}\right)a_{11}b_{21} & a_{12}b_{22}a_{11}b_{21}\\
\displaystyle
0 & 0 \\
\displaystyle
0  & 0 \\
\end{array}
\right).
\end{eqnarray*}
Now because of (\ref{3}) and by applying the elementary row transformation, new(row 2)=row 2 - row 1,
our $D$ matrix reduces to
\begin{eqnarray*}
\displaystyle
D=\left(
\begin{array}{cc}
\displaystyle
-b_{11}\left(a_{21}\right)a_{12}b_{22} & a_{11}b_{21}a_{12}b_{22}\\
\displaystyle
0   & 0\\
\displaystyle
0 & 0 \\
\displaystyle
0 & 0 \\
\end{array}
\right).
\end{eqnarray*}
Hence, rank$(D) = 1$.
This completes our proof.
$\square$

Next we show that rank$(D)$ is bigger than one if $A$ and $B$ are not compatible.

\noindent
{\bf Corollary 2:}
In a $(2\times 2)$ case, rank$(D) = 2$, if $A$ and $B$ are not compatible.

\noindent
{\bf Proof:}
This actually follows from our previous result.
After some elementary row-transformations, our  $D$ matrix reduces to
\begin{eqnarray*}
\displaystyle
D=\left(
\begin{array}{cc}
\displaystyle
-b_{11}\left(a_{21}\right)a_{12}b_{22} & a_{11}b_{21}a_{12}b_{22}\\
\displaystyle
-b_{12}\left(a_{22}\right)a_{11}b_{21} & a_{12}b_{22}a_{11}b_{21}\\
\displaystyle
0  & 0 \\
\displaystyle
0  & 0 \\
\end{array}
\right).
\end{eqnarray*}
When $A$ and $B$ are incompatible, (\ref{3}) does not hold so the rows of $D$ are not proportional.
Therefore, $\text{\rm rank}(D)=2$.
This completes our proof.
$\square$

\subsection{Examples of $(2\times 2)$  case}

First of   all we consider a situation where the
two matrices  $A$ and  $B$ are compatible with all the elements strictly positive.
\begin{itemize}

\item
Suppose we have
\begin{eqnarray*}
\displaystyle
A=\left(
\begin{array}{cc}
1/4&2/3\\
3/4&1/3\\
\end{array}
\right)
\end{eqnarray*}
and
\begin{eqnarray*}
\displaystyle
B=\left(
\begin{array}{cc}
1/3  & 2/3 \\
3/4 & 1/4 \\
\end{array}
\right).
\end{eqnarray*}
In this case, the resulting $D$ matrix is
\begin{eqnarray*}
\displaystyle
D=\left(
\begin{array}{cc}
   -0.2500000&  0.1875000\\
 -0.2222222 & 0.1666667\\
  0.2500000 &-0.1875000\\
  0.2222222 &-0.1666667\\
\end{array}
\right).
\end{eqnarray*}
In this case, rank$(D) = 1$.
So, $A$ and $B$ are compatible as can be
verified by checking the cross product ratios of $A$ and $B$.
The solution for the joint probability distribution in this case is
\begin{eqnarray*}
\displaystyle
P=\left(
\begin{array}{cc}
1/7  & 2/7 \\
3/7 & 1/7 \\
\end{array}
\right).
\end{eqnarray*}

\item
Next, we consider two matrices $A$ and $B$ of the following forms:
\begin{eqnarray*}
\displaystyle
A=\left(
\begin{array}{cc}
1/7&3/4\\
6/7&1/4\\
\end{array}
\right)
\end{eqnarray*}
and
\begin{eqnarray*}
\displaystyle
B=\left(
\begin{array}{cc}
2/5  & 3/5 \\
3/8 & 5/8 \\
\end{array}
\right).
\end{eqnarray*}
The resulting $D$ matrix in this case is
\begin{eqnarray*}
\displaystyle
D=\left(
\begin{array}{cc}
 -0.3428571 & 0.05357143\\
-0.1500000  &0.46875000\\
  0.3428571 &-0.05357143\\
  0.1500000 &-0.46875000\\
\end{array}
\right).
\end{eqnarray*}
In this case, rank$(D) = 2$, so $A$ and $B$ are incompatible.

\end{itemize}

\subsection{Compatibility in $(3\times 3)$ case}

\begin{enumerate}

\item
First of all let us consider a compatible case of type 1,
where by type 1, we mean all the elements of $A$ and $B$ are strictly positive.

Let
\begin{eqnarray*}
\displaystyle
A=\left(
\begin{array}{ccc}
  1/5& 2/7& 3/8\\
  3/5 &2/7& 1/8\\
  1/5& 3/7& 1/2\\
\end{array}
\right)
\end{eqnarray*}
and
\begin{eqnarray*}
\displaystyle
B=\left(
\begin{array}{ccc}
  1/6& 1/3& 1/2\\
  1/2 &1/3& 1/6\\
  1/8& 3/8& 1/2\\
\end{array}
\right).
\end{eqnarray*}
In this case, the corresponding $D$ matrix is
\begin{eqnarray*}
\displaystyle
D^{(9\times 3)}=\left(
\begin{array}{ccccccccc}
  -0.13333333 & 0.10000000 & 0.0250000\\
  -0.23809524 & 0.09523810 & 0.1071429\\
  -0.31250000  &0.06250000  &0.1875000\\
   0.10000000 &-0.20000000  &0.0750000\\
  0.09523810 &-0.23809524  &0.1071429\\
   0.06250000 &-0.14583333  &0.0625000\\
   0.03333333  &0.10000000  &-0.1000000\\
  0.14285714  &0.14285714   &-0.2142857\\
   0.25000000 & 0.08333333 &-0.2500000\\
\end{array}
\right).
\end{eqnarray*}
Note that in this case rank$(D) = 2$ and hence  $A$ and  $B$ are compatible.
Now solving the  equation as mentioned earlier
the solution for the marginal of  $X$ is $\underline{\eta}=(0.3, 0.3, 0.4)$.

\item
Next consider a compatible case of type 2, where by type 2,
we mean some of the elements of $A$ and $B$
are zeros which appear in the same positions in both $A$ and $B$.
Suppose $A$ and $B$ take the following forms:
\begin{eqnarray*}
\displaystyle
A=\left(
\begin{array}{ccc}
  1/3& 0& 2/3\\
  0 &1/2& 1/3\\
  2/3& 1/2& 0\\
\end{array}
\right)
\end{eqnarray*}
and
\begin{eqnarray*}
\displaystyle
B=\left(
\begin{array}{ccc}
  1/3& 0& 2/3\\
  0 &1/2& 1/2\\
  2/3& 1/3& 0\\
\end{array}
\right).
\end{eqnarray*}
In this case, the corresponding  $D$ matrix is
\begin{eqnarray*}
\displaystyle
D^{(9 \times 3)}=\left(
\begin{array}{ccc}
-0.2222222 &  0.0000000 &  0.2222222\\
  0.0000000 & 0.0000000 & 0.0000000\\
  -0.2222222&  0.3333333&  0.0000000\\
   0.0000000&  0.0000000&  0.0000000\\
   0.0000000& -0.2500000 & 0.1666667\\
  0.2222222 & -0.3333333 &  0.0000000\\
   0.2222222&  0.0000000 & -0.2222222\\
   0.0000000&  0.2500000 & -0.1666667\\
   0.0000000&  0.0000000 &  0.0000000\\
\end{array}
\right).
\end{eqnarray*}
In this case, rank$(D) = 2$ and indeed $A$ and $B$ are compatible.
Solving the  equation, we obtain the solution for the marginal
of $X$ as $\underline{\eta}=(0.375, 0.250, 0.375)$.
The joint probability distribution in this case is
\begin{eqnarray*}
\displaystyle
P=\left(
\begin{array}{ccc}
  0.125& 0.000& 0.250\\
 0.000 &0.125& 0.125\\
 0.250& 0.125& 0.000\\
\end{array}
\right).
\end{eqnarray*}

\item
Next we consider an incompatible case of type 1,
in which $A$ and $B$ have the following forms:
\begin{eqnarray*}
\displaystyle
A=\left(
\begin{array}{ccc}
  0.2&  0.3&  0.1\\
  0.1 & 0.4 & 0.4\\
  0.7 & 0.3 & 0.5\\
\end{array}
\right)
\end{eqnarray*}
and
\begin{eqnarray*}
\displaystyle
B=\left(
\begin{array}{ccc}
  0.2&  0.1&  0.7\\
  0.3 & 0.4 & 0.3\\
  0.1 & 0.4 & 0.5\\
\end{array}
\right).
\end{eqnarray*}
In this case, the resulting $D$ matrix is
\begin{eqnarray*}
\displaystyle
D^{(9 \times 3)}=\left(
\begin{array}{ccc}
 -0.16 & 0.06 & 0.02\\
  -0.07&  0.12 & 0.12\\
  -0.63 & 0.03 & 0.05\\
   0.02 &-0.27 & 0.01\\
   0.04 &-0.24  &0.16\\
   0.28 &-0.18 & 0.20\\
   0.14 & 0.21& -0.03\\
   0.03 & 0.12& -0.28\\
   0.35 & 0.15& -0.25\\
\end{array}
\right).
\end{eqnarray*}
In this case,  $\text{\rm rank}(D)=3$  and hence  $A$  and  $B$  are not compatible.

\item
Next, we consider an incompatible case of type 2 in which $A$ and $B$ have the following forms:
\begin{eqnarray*}
\displaystyle
A=\left(
\begin{array}{ccc}
  0&  1/3&  0\\
  1 & 1/3 & 1/2\\
  0 & 1/3 & 1/2\\
\end{array}
\right)
\end{eqnarray*}
and
\begin{eqnarray*}
\displaystyle
B=\left(
\begin{array}{ccc}
  0&  1&  0\\
  1/4 & 1/2 & 1/4\\
  0 & 1/5 & 4/5\\
\end{array}
\right).
\end{eqnarray*}
In this case, the resulting $D$ matrix is
\begin{eqnarray*}
\displaystyle
D^{(9 \times 3)}=\left(
\begin{array}{ccc}
 0.0000000 & 0.0000000&  0.00000000\\
  -0.6666667 & 0.1666667 & 0.06666667\\
   0.0000000 & 0.0000000 & 0.00000000\\
   0.0000000  &0.0000000 & 0.00000000\\
  0.3333333 &-0.3333333 & 0.06666667\\
   0.0000000 &-0.1250000 & 0.40000000\\
  0.0000000  &0.0000000 & 0.00000000\\
   0.3333333 & 0.1666667 &-0.13333333\\
   0.0000000 & 0.1250000 &-0.40000000\\
\end{array}
\right).
\end{eqnarray*}
In this case,  $\text{\rm rank}(D)=3$  and hence  $A$  and  $B$  are not compatible.
\end{enumerate}

\subsection{Proof that rank$(D)=2$ when $A$ and $B$ are compatible in a  $(3\times 3)$ case}

\begin{theorem}
For a $3\times 3$ case, if the two matrices $A$ and $B$ are compatible then $\text{\rm rank(D)}=I-1=2$.
\end{theorem}

\noindent
{\bf Proof:}
The form of the $D$ matrix in a $(3\times 3)$ case is
\begin{eqnarray*}
\displaystyle
D^{(9 \times 3)}=\left(
\begin{array}{ccc}
\displaystyle
b_{11}\left(a_{11}-1\right)&  a_{11}b_{21}& a_{11}b_{31}\\
\displaystyle
b_{12}\left(a_{12}-1\right)& a_{12}b_{22}& a_{12}b_{32}\\
\displaystyle
b_{13}\left(a_{13}-1\right)& a_{13}b_{23}& a_{13}b_{33}\\
\displaystyle
a_{21}b_{11} & b_{21}\left(a_{21}-1\right) & a_{21}b_{31} \\
\displaystyle
a_{22}b_{12}& b_{22}\left(a_{22}-1\right)& a_{22}b_{32}\\
\displaystyle
a_{23}b_{13}& b_{23}\left(a_{23}-1\right) & a_{23}b_{33}\\
\displaystyle
a_{31}b_{11} & a_{31}b_{21} & b_{31}\left(a_{31}-1\right)\\
\displaystyle
a_{32}b_{12} & a_{32}b_{22} & b_{32}\left(a_{32}-1\right)\\
\displaystyle
a_{33}b_{13} &a_{33}b_{23}& b_{33}\left(a_{33}-1\right) \\
\end{array}
\right).
\end{eqnarray*}
However, if matrices $A$ and $B$ are compatible then all possible cross
product ratios of $A$  are equal to the corresponding cross product ratios of $B$.
First of all we apply the following elementary row operations:
\begin{itemize}

\item
new(row 1)=row 1 + row 4 + row 7;

\item
new(row 2)=row 2 + row 5 + row 8;

\item
new(row 3)=row 3 + row 6 + row 9.

\end{itemize}
So, the matrix $D$ reduces to
\begin{eqnarray*}
\displaystyle
D=\left(
\begin{array}{ccccccccc}
\displaystyle
0&0&0\\
\displaystyle
0&0&0\\
\displaystyle
0&0&0\\
\displaystyle
a_{21}b_{11} &b_{21}\left(a_{21}-1\right) &a_{21}b_{31} \\
\displaystyle
a_{22}b_{12}&b_{22}\left(a_{22}-1\right)&a_{22}b_{32}\\
\displaystyle
a_{23}b_{13}&b_{23}\left(a_{23}-1\right) &a_{23}b_{33}\\
\displaystyle
a_{31}b_{11} &a_{31}b_{21} &b_{31}\left(a_{31}-1\right)\\
\displaystyle
a_{32}b_{12} &a_{32}b_{22} & b_{32}\left(a_{32}-1\right)\\
\displaystyle
a_{33}b_{13} &a_{33}b_{23}&b_{33}\left(a_{33}-1\right) \\
\end{array}
\right).
\end{eqnarray*}
Next consider the following elementary row and column operations:
\begin{itemize}

\item
new(row 4)=$\frac {\text{\rm row} 4}{a_{21}}$;

\item
new(row 5)=$\frac {\text{\rm row} 5}{a_{22}}$;

\item
new(row 6)=$\frac {\text{\rm row} 4}{a_{23}}$;

\item
new(row 7)=$\frac {\text{\rm row} 7}{a_{31}}$;

\item
new(row 8)=$\frac {\text{\rm row} 8}{a_{32}}$;

\item
new(row 9)=$\frac {\text{\rm row} 4}{a_{33}}$;

\item
new(col 4)=col 4 + col 5 + col 6;

\item
new(col 7)=col 7 + col 8 + col 9.

\end{itemize}
So, the $D$ matrix is:
\begin{eqnarray*}
\displaystyle
D^{(9 \times 3)}=\left(
\begin{array}{ccc}
\displaystyle
0&0&0\\
\displaystyle
0&0&0\\
\displaystyle
0&0&0\\
\displaystyle
1 &1-\left(\frac {b_{21}}{a_{21}}+\frac {b_{22}}{a_{22}}+\frac {b_{23}}{a_{23}}\right) &1 \\
\displaystyle
b_{12}&b_{22}\left(1-\frac {1}{a_{22}}\right)&b_{32}\\
\displaystyle
b_{13}&b_{23}\left(1-\frac {1}{a_{23}}\right) &b_{33}\\
\displaystyle
1 &1 &1-\left(\frac {b_{31}}{a_{31}}+\frac {b_{32}}{a_{32}}+\frac {b_{33}}{a_{33}}\right)\\
\displaystyle
b_{12} &b_{22} & b_{32}\left(1-\frac {1}{a_{32}}\right)\\
\displaystyle
b_{13} &b_{23}&b_{33}\left(1-\frac {1}{a_{33}}\right) \\
\end{array}
\right).
\end{eqnarray*}
Then we consider the following:
\begin{itemize}

\item
new(row 5)=row 5 + row 6;

\item
new(row 8)=row 8 + row 9;

\item
new(row 5)=row 5 - row 8;

\item
new(row 4)=row 4 - row 7;

\item
new(row 6)=row 6 - row 8.

\end{itemize}
Then our $D$ matrix reduces to
\begin{eqnarray*}
\displaystyle
D^{(9 \times 3)}=\left(
\begin{array}{ccc}
\displaystyle
0&0&0\\
\displaystyle
0&0&0\\
\displaystyle
0&0&0\\
\displaystyle
0 &-\left(\frac {b_{21}}{a_{21}}+\frac {b_{22}}{a_{22}}+\frac {b_{23}}{a_{23}}\right)
\displaystyle
&\left(\frac {b_{31}}{a_{31}}+\frac {b_{32}}{a_{32}}+\frac {b_{33}}{a_{33}}\right) \\
\displaystyle
0&-\left(\frac {b_{22}}{a_{22}}+\frac {b_{23}}{a_{23}}\right)&\left(\frac {b_{32}}{a_{32}}+\frac {b_{33}}{a_{33}}\right)\\
\displaystyle
0&-\frac {b_{23}}{a_{23}} &\frac {b_{33}}{a_{33}}\\
\displaystyle
1 &1 &1-\left(\frac {b_{31}}{a_{31}}+\frac {b_{32}}{a_{32}}+\frac {b_{33}}{a_{33}}\right)\\
\displaystyle
-b_{11} &-b_{21} &-b_{31}-\left(\frac {b_{32}}{a_{32}}+\frac {b_{33}}{a_{33}}\right) \\
\displaystyle
b_{13} &b_{23}&b_{33}\left(1-\frac {1}{a_{33}}\right) \\
\end{array}
\right).
\end{eqnarray*}
Again we consider new (row 4)=row 4 - row 5,  new(row 5)=row 5 - row 6,  so the $D$ matrix reduces to
\begin{eqnarray*}
\displaystyle
D^{(9 \times 3)}=\left(
\begin{array}{ccccccccc}
\displaystyle
0&0&0\\
\displaystyle
0&0&0\\
\displaystyle
0&0&0\\
\displaystyle
0 &-\frac {b_{21}}{a_{21}} &\frac {b_{31}}{a_{31}} \\
\displaystyle
0&-\frac {b_{22}}{a_{22}}&\frac {b_{32}}{a_{32}}\\
\displaystyle
0&-\frac {b_{23}}{a_{23}} &\frac {b_{33}}{a_{33}}\\
\displaystyle
1 &1 &1-\left(\frac {b_{31}}{a_{31}}+\frac {b_{32}}{a_{32}}+\frac {b_{33}}{a_{33}}\right)\\
\displaystyle
-b_{11} &-b_{21} &-b_{31}-\left(\frac {b_{32}}{a_{32}}+\frac {b_{33}}{a_{33}}\right) \\
\displaystyle
b_{13} &b_{23}&b_{33}\left(1-\frac {1}{a_{33}}\right) \\
\end{array}
\right).
\end{eqnarray*}
Note that in this case we have rank$\left(D^{(9 \times 3)}\right)$ $\leq \text{\rm min}(9,3)=3$.
However, for our matrix $D$, the determinants of all possible submatrices of order $(3\times 3)$
are zero and hence rank$\left( D^{(9 \times 3)} \right) < 3$.
We know that the rank of a matrix is the
highest order non-vanishing determinant.
Let us consider the determinant of any submatrix of order $(2\times 2)$:
\begin{eqnarray*}
\displaystyle
B=\left(
\begin{array}{cc}
\displaystyle
-b_{11} &-b_{21}\\
\displaystyle
b_{13} &b_{23}\\
\end{array}
\right).
\end{eqnarray*}
The determinant for $B$ is
\begin{eqnarray*}
\displaystyle
\text{\rm det}(B) = -b_{11}b_{23} + b_{13}b_{21} \neq 0.
\end{eqnarray*}
So, we have  $\text{\rm rank}(D)=I-1=2$, which follows from the definition of the rank of a matrix.
Hence, $A$ and $B$ are compatible if and only if rank$(D)=I-1$.
However, if $A$ and $B$ are not compatible then the rows of $A$ are not proportional
to the rows of $B$ which implies that $\text{\rm rank}(D)>2$.
Hence, the proof.
$\square$

Summarizing Sections 5.1 to 5.4, we can say the following:
i) if $A$ and $B$ both $2 \times 2$ are compatible then the rank of $D$ is $1$;
ii) if $A$ and $B$ both $2 \times 2$ are not compatible then the rank of $D$ is $2$;
iii) if $A$ and $B$ both $3 \times 3$ are compatible then the rank of $D$ is $2$.

\section{Application to a real life data}

Here, we illustrate our proposed method by considering a real life example.
The data has been collected from a local grocery store in Clarksville, Tennessee regarding the
supply ($X$) and demand ($Y$) of a particular brand of cigarette in big cartons in a particular week.
We denote the possible values of $X$ as $X=1,2,3$ and the possible values of $Y$ as $Y=0,1,2$.
The conditional probabilities are:
\begin{eqnarray*}
\displaystyle
A=\left(
\begin{array}{ccc}
  0.2& 0.3& 0.1\\
  0.1&0.4& 0.4\\
 0.7& 0.3& 0.5\\
\end{array}
\right)
\end{eqnarray*}
and
\begin{eqnarray*}
\displaystyle
B=\left(
\begin{array}{ccc}
  0.2& 0.1& 0.7\\
  0.3&0.4& 0.3\\
 0.1& 0.4& 0.5\\
\end{array}
\right).
\end{eqnarray*}
In this case, the resulting $D$ matrix is
\begin{eqnarray*}
\displaystyle
D^{(9 \times 3)}=\left(
\begin{array}{ccc}
 -0.16 & 0.06 & 0.02\\
  -0.07&  0.12 & 0.12\\
  -0.63 & 0.03 & 0.05\\
   0.02 &-0.27 & 0.01\\
   0.04 &-0.24  &0.16\\
   0.28 &-0.18 & 0.20\\
   0.14 & 0.21& -0.03\\
   0.03 & 0.12& -0.28\\
   0.35 & 0.15& -0.25\\
\end{array}
\right).
\end{eqnarray*}
In this case,  $\text{\rm rank}(D)=3$  and hence  $A$  and  $B$  are not compatible.
So, in this situation it is not possible to obtain a perfect joint distribution of  $X$ and $Y$.
However, one may obtain a most nearly compatible joint distribution $P$, see Ghosh and Balakrishnan (2013).

\section{Concluding remarks}

Compatible conditional specifications of distributions are of
fundamental importance in modeling scenarios.
Moreover, in Bayesian prior elicitation contexts, inconsistent conditional specifications are to be expected.
In such situations, interest  centers on most nearly compatible distributions.
In the finite discrete case, a variety of compatibility conditions can be derived.
In this paper, we have discussed in detail the problem of compatibility
by identifying it as a LP problem and we have developed a rank based criterion.
We have shown that the rank of the matrix  (whose elements are constructed from the two given matrices $A$ and $B$)
under compatibility is $I-1$ for a ($2\times 2$).
A similar result holds for any dimension $(I\times J)$.
A significant amount of the material for this paper draws heavily on
Arnold et al. (1999) and Arnold and Gokhale (1998).

A future work is to check if the results in this paper can
provide an ideal to construct an approximated joint distribution
that is ``logically acceptable''
when the given conditional distributions are not compatible.
Furthermore, if three or more probability distributions are pairwise compatible,
what conditions do we need for collective compatibility in terms of our proposed LP approach?
Another future work is to give statistical insights to the results.

\section*{Acknowledgments}

The authors would like to thank the Editor and the referee
for careful reading and for their comments which greatly improved the paper.


\begin{thebibliography}{999}


\bibitem{}
Amemiya, T. (1975).
Qualitative response models.
Annals of Economic Social Measurement, 4, 363-372.




\bibitem{}
Arnold, B.C., Castillo, E.,  and Sarabia, J.M. (1999).
Conditional Specification of Statistical Models.
Springer Verlag, New York.


\bibitem{}
Arnold, B.C., Castillo, E., and  Sarabia, J.M. (2001).
Quantification of incompatibility of conditional and marginal information.
Communications in Statistics---Theory and Methods, 30, 381-395.


\bibitem{}
Arnold, B.C., Castillo, E., and  Sarabia, J.M. (2002).
Exact and near compatibility of discrete conditional distributions.
Computational Statistics and Data Analysis, 40, 231-252.


\bibitem{}
Arnold, B.C., Castillo, E., and  Sarabia, J.M. (2004).
Compatibility of partial or complete conditional probability specifications.
Journal of Statistical Planning and Inference, 123, 133-159.



\bibitem{}
Arnold, B.C. and Gokhale, D.V. (1994).
On uniform marginal representations of contingency tables.
Statistics and Probability Letters, 21, 311-316.


\bibitem{}
Arnold, B.C. and Gokhale, D.V. (1998).
Distributions most nearly compatible with given families of conditional distributions.
Test, 7, 377-390.


\bibitem{}
Arnold, B.C. and Press, S.J. (1989).
Compatible conditional distributions.
Journal of the American Statistical Association, 84, 152-156.


\bibitem{}
Cacoullos, T. and Papageorgiou, H. (1983).
Characterizations of discrete distributions by a conditional distribution and a regression function.
Annals of the Institute of Statistical Mathematics, 35, 95-103.





\bibitem{}
Ghosh, I. and Balakrishnan, N. (2013).
Study of incompatibility or near compatibility of bivariate discrete conditional probability distributions through divergence measures.
Journal of Statistical Computation and Simulation, 84, 1-14.


\bibitem{}
Gourieroux, C. and Montfort, A. (1979).
On the characterization of a joint probability distribution by conditional distributions.
Journal of  Econometrics, 10, 115-118.



\bibitem{}
Ip, E.H. and Wang, Y.J. (2009).
Canonical representation of conditionally specified multivariate discrete distributions.
Journal of Multivariate Analysis, 100, 1282-1290.


\bibitem{}
Kuo, K.-L. and Wang, Y.J. (2011).
A simple algorithm for checking compatibility among discrete conditional distributions.
Computational Statistics and Data Analysis, 55, 2457-2462.


\bibitem{}
Lauritzen, S.L. and Spiegelhalter, D.J. (1988).
Local computations with probabilities on
graphical structures and their application to expert systems.
Journal of the Royal Statistical Society, B, 50, 157-194.


\bibitem{}
Nerlove, M. and Press, S.J. (1986).
Multivariate log-linear probability models in econometrics.
In: Advances in Statistical Analysis and Statistical Computing,
editor R.S. Mariano, pp. 117-171, JAI Press, Greenwich, Connecticut.


\bibitem{}
Patil, G.P. (1965).
Certain characteristic properties of multivariate discrete probability distributions
akin to the Bates-Neyman model in the theory of accident proneness.
Sankhy\=a, A, 27, 259-270.




\bibitem{}
Pearl, J. (1988).
Probabilistic Reasoning in Intelligent Systems: Networks of Plausible Inference.
Morgan Kaufmann, San Mateo, California.



\bibitem{}
Song, C.-C., Li, L.-A., Chen, C.-H., Jiang, T.J. and Kuo, K.-L. (2010).
Compatibility of finite discrete conditional distributions.
Statistica Sinica, 20, 423-440.


\bibitem{}
Tian, G.-L., Tan, M., Ng, K.W. and Tang, M.-L. (2009).
A unified method for checking compatibility and uniqueness for finite discrete conditional distributions.
Communications in Statistics---Theory and Methods, 38, 115-129.


\bibitem{}
Wang, Y.J. (2012).
Comparisons of three approaches for discrete conditional models.
Communications in Statistics---Simulation and Computation, 41, 32-43.



\bibitem{}
Wang, Y.J. and Kuo, K.-L. (2010).
Compatibility of discrete conditional distributions with structural zeros.
Journal of Multivariate Analysis, 101, 191-199.



\bibitem{}
Yao, Y.-C., Chen, S.-C. and Wang, S.-H. (2014).
On compatibility of discrete full conditional distributions: A graphical representation approach.
Journal of Multivariate Analysis, 124, 1-9.



\end{thebibliography}
\end{document}